\documentclass[11pt]{amsart}
\usepackage[letterpaper,margin=1in]{geometry}                
\usepackage{amssymb,amsmath,amsfonts,stmaryrd,mathscinet}
\usepackage{xcolor}
\usepackage[pdfproducer={pdfLaTeX},
            pdfpagemode=None,
            bookmarksopen=true
            bookmarksnumbered=true]{hyperref}
\usepackage{graphicx}
\usepackage{multicol}

\usepackage{tikz}
\usetikzlibrary{decorations.markings,arrows,decorations.pathreplacing}
\usepackage{arydshln}

\allowdisplaybreaks

\newcommand\arXiv[1]{\href{http://arxiv.org/abs/#1}{\nolinkurl{arXiv:#1}}}
\newcommand\MRnumber[1]{\href{http://www.ams.org/mathscinet-getitem?mr=#1}{\nolinkurl{MR#1}}}
\newcommand\DOI[1]{\href{http://dx.doi.org/#1}{\nolinkurl{DOI:#1}}}
\newcommand\MAILTO[1]{\href{mailto:#1}{\nolinkurl{#1}}}

\newcommand\mytheorem[2][Theorem]{
\vspace*{8pt}\noindent\textbf{#1.} \emph{#2}

\vspace*{8pt}}

\usepackage{dsfont}
\renewcommand\mathbb\mathds

\newcommand\bC{\mathbb C}

\newcommand\bF{\mathbb F}

\newcommand\bK{\mathbb K}
\newcommand\bL{\mathbb L}
\newcommand\bM{\mathbb M}
\newcommand\bN{\mathbb N}

\newcommand\bP{\mathbb P}
\newcommand\bQ{\mathbb Q}
\newcommand\bR{\mathbb R}

\newcommand\bZ{\mathbb Z}

\newcommand\cC{\mathcal C}

\newcommand\cO{\mathcal O}

\newcommand\cZ{\mathcal Z}

\newcommand\rK{\mathrm K}

\usepackage[mathscr]{euscript}

\DeclareMathOperator\homology{H}
\renewcommand\H{\homology}

\newcommand\bmu{\boldsymbol{\mu}}

\newcommand\longto\longrightarrow
\newcommand\mono\hookrightarrow
\newcommand\epi\twoheadrightarrow
\newcommand\isom{\overset\sim\to}
\newcommand\<\langle
\renewcommand\>\rangle
\newcommand\sminus\smallsetminus

\newcommand\st{\text{ s.t.\ }}
\newcommand\op{\mathrm{op}}
\newcommand\id{\mathrm{id}}

\DeclareMathOperator\Aut{Aut}

\DeclareMathOperator\Gal{Gal}
\DeclareMathOperator\Out{Out}

\DeclareMathOperator\Br{Br}

\newcommand\Vect{\cat{Vect}}
\newcommand\Rep{\cat{Rep}}

\newcommand\define[1]{\emph{#1}}
\newcommand\cat[1]{\textsc{#1}}

\title{Galois action on VOA gauge anomalies}
\author{Theo Johnson-Freyd}
\email{\MAILTO{theojf@pitp.ca}}
\address{Perimeter Institute for Theoretical Physics, Waterloo, Ontario}
\thanks{This project was started during the CIRM conference ``Representation Theory, Mathematical Physics and Integrable Systems'' celebrating the 60th birthday of Kolya Reshetikhin. It was from Kolya that I first learned about  (quantum) Hamiltonian reduction, spin chains, conformal field theory, and modular tensor categories, and I thank him for his ongoing interest and inspiration. I also thank the anonymous referee for their thoughtful comments and corrections.
Research at the Perimeter Institute for Theoretical Physics is supported by the Government of Canada through the Department of Innovation, Science and Economic Development Canada and by the Province of Ontario through the Ministry of Research, Innovation and Science.}

\begin{document}

\maketitle

\begin{center}
\emph{For Kolya Reshetikhin on the occasion of his 60th birthday.}
\end{center}

\begin{abstract}
  Assuming regularity of the fixed subalgebra, any action of a finite group $G$ on a holomorphic VOA $V$ determines a gauge anomaly  $\alpha \in \H^3(G; \bmu)$, where $\bmu \subset \bC^\times$ is the group of roots of unity. We show that under Galois conjugation $V \mapsto {^\gamma V}$, the gauge anomaly transforms as $\alpha \mapsto \gamma^2(\alpha)$. This provides an {a priori} upper bound of $24$ on the order of anomalies of actions preserving a $\bQ$-structure, for example the Monster group $\bM$ acting on its Moonshine VOA $V^\natural$. We speculate that each field $\bK$ should have a ``vertex Brauer group'' isomorphic to $\H^3(\Gal(\bar{\bK}/\bK); \bmu^{\otimes 2})$. In order to motivate our constructions and speculations, we warm up with a discussion of the ordinary Brauer group, emphasizing the analogy between VOA gauging and quantum Hamiltonian reduction.
\end{abstract}

\vspace*{18pt}
This article studies the behavior of gauge (aka 't Hooft) anomalies under the action of Galois automorphisms. Our  Main Theorem, stated and proved in \S\ref{sec:mainthm}, asserts that if $V$ is a holomorphic VOA over $\bC$ on which a finite group $G$ acts with anomaly $\alpha$, and if $\gamma$ is a field automorphism of~$\bC$, then $G$ acts on the Galois-conjugate VOA $^\gamma V$ with anomaly $\gamma^2(\alpha)$, rather than, as might be naively expected, $\gamma(\alpha)$. One corollary is \cite[Conjecture~1]{MR3916985}: there is an {a priori} upper bound of $24$ on the order of gauge anomalies of actions preserving a $\bQ$-form of the VOA.  In \S\ref{sec:conclusions} we list some further questions and conjectures, including the possibility of a ``vertex Brauer group.''

These results are conditional on a widely believed Regularity Assumption, which asserts that 
the representation theory of the $G$-fixed sub-VOA $V^G$ is finite semisimple.
The Regularity Assumption is required in order for gauge anomalies to be well-defined. We review the construction of gauge anomalies for VOAs in~\S\ref{sec:VOAs}.
As a warm-up,~\S\ref{sec:warmup} discusses gauge anomalies in the case of actions on Azumaya algebras: anomalies live in $\H^2(G;\bK^\times)$, with the standard Galois action; such anomalies lead quickly to the  Galois-cohomological description of the Brauer group as $\Br(\bK) \cong \H^2(\Gal(\bar{\bK}/\bK);\bar{\bK}^\times)$.

\section{Gauge anomalies for associative algebras}\label{sec:warmup}

The mathematics in this section is not new. Its goal is to present a familiar story in a way that will generalize to VOAs. A good source for the complete story is \cite{MR2266528}.

\subsection{Azumaya algebras}
A $\bK$-algebra $A$ is \define{Azumaya} if it is finite-dimensional and simple, and its centre is exactly $\bK$. Let $A^\op$ denote the opposite algebra. Then $A$ is Azumaya if and only if its \define{enveloping algebra} $A \otimes A^\op$ is isomorphic to a matrix algebra. This is in turn equivalent to $A \otimes A^\op$ being Morita-trivial. Modules for $A \otimes A^\op$ are precisely $A$-$A$-bimodules, and so an Azumaya algebra is an algebra whose (monoidal) category of bimodules is  equivalent to $\Vect$.

\subsection{Physical interpretation}\label{subsec:azumaya-physics}
Azumaya algebras have the following interpretation in quantum physics. Consider a (one-dimensional) chain of atoms. To each atom we associate a ``local Hilbert space'' given by the underlying vector space of~$A$ (so that a basis of $A$ are the ``spins'' of the ``spin chain'');
  if the atoms are index by $i \in \{\dots, 1,2,\dots\}$, then we will write $A_i$ for the copy of $A$ at location $i$.
 The total Hilbert space is the tensor product $\bigotimes_i A_i$ of all the local Hilbert spaces. (That tensor product makes sense only for chains of finite length.) 
 
 We now place a Hamiltonian on this Hilbert space, which will be of ``commuting projector'' type. Specifically, consider the multiplication map $m : A \otimes A \to A$. It is a map of $A$-$A$-bimodules, and by our Azumaya condition corresponds to a linear map $B \to \bK$, where we have used the identification  $\cat{Bimod}(A) \simeq \Vect$ to identify ${_A A_A} \in \cat{Bimod}(A)$ with $\bK \in \Vect$ and to identify ${_A A \otimes A_A} \in \cat{Bimod}(A)$ with some vector space $B \in \Vect$. (The vector space $B$ will be noncanonically isomorphic to the underlying vector space of $A$.) The map $B \to \bK$ is nonzero, hence a surjection, and so can be split: we can choose a bimodule map $\Delta : {_A A_A} \to {_A A \otimes A_A}$ such that $m \circ \Delta = \id_A$. Then $\Delta \circ m : A^{\otimes 2} \to A^{\otimes 2}$ is a projection (it squares to itself) onto a vector space isomorphic to $A$. 
 Let $\id_i$ denote the identity on $A_i$, and for a set $I$ of indices let $\id_I = \bigotimes_{i\in I} \id_i$. Also write $(\Delta \circ m)_{i,i+1}$ for the map $\Delta \circ m$ acting on $A_i \otimes A_{i+1}$. 
 The Hamiltonian we will use is:
 $$\sum_{i} \id_{<i} \otimes (\id - \Delta \circ m)_{i,i+1} \otimes \id_{>i+1}.$$

 This Hamiltonian is gapped (i.e.\ it has a finite-dimensional lowest eigenspace, and its spectrum is discrete), and its vector space of ground states  is isomorphic to $A$. Since the chain of atoms was $1$-space-dimensional, together with the Hamiltonian we have a ``gapped topological system of $(1+1)$-dimensional matter.'' 
 It is expected that gapped physical systems admit low-energy effective descriptions in terms of TQFTs. For the spin-chain model built from $A$, the low-energy TQFT is the one constructed from $A$ a la \cite{Schommer-Pries:thesis}.

The gapped topological phase (meaning the connected component of the space of all gapped systems under deformations that do not close the gap) constructed from a semisimple algebra $A$ depends only on the Morita equivalence class of $A$
(see e.g.\ \cite[\S2.4]{GJFa}).
To encode the algebra $A$ itself, one can remember not just the phase of $(1+1)$-dimensional matter, but also a choice of boundary condition. Namely, the spin-chain constructed above has an obvious ``free'' boundary condition, whose
boundary observables form a copy of the associative algebra $A$. The bulk observables live in the centre $\cZ(A)$, which is trivial by assumption. 
Thus the bulk theory is almost trivial:
it is ``short-range entangled.''
The class of $A$ in the Brauer group $\Br(\bK)$ can be considered a ``gravitational anomaly'' for the boundary system: the boundary system fails to be an absolute $(0+1)$-dimensional QFT only because of the Brauer class of $A$.
 Some discussion of absolute versus anomalous quantum systems, in the framework of TQFTs, is available in \cite{MR3165462,me-Heisenberg}; for a condensed matter discussion focussing on gauge anomalies, see \cite[\S7]{GJFa}.

\subsection{Construction of gauge anomalies}\label{assoc-construction}
Let $A$ be Azumaya and $G \subset \Aut(A)$ a finite group of automorphisms of $A$. (One can be more general, but we will not.) We will associate to these data a \define{gauge anomaly} $\alpha \in \H^2(G; \bK^\times)$, where $G$ acts trivially on $\bK^\times$.

There are many ways to define $\alpha$. The quickest is to consider the exact sequence, valid for any algebra,
\begin{equation}
 \{1\} \to \cZ(A^\times) \to A^\times \to \Aut(A) \to \Out(A) \to \{1\} , \label{eqn-assoc}
\end{equation}
where the map $A^\times \to \Aut(A)$ picks out the inner automorphisms, and $\Out(A)$ denotes the group of outer automorphisms of $A$. Two simplifications occur when $A$ is Azumaya. First, from the definition in terms of central simple algebras, we see that $\cZ(A^\times) = \bK^\times$. Second, $\Out(A)$ is trivial by the Skolem--Noether theorem. Let us quickly prove the second claim. Any automorphism $\phi : A\to A$ determines an invertible bimodule $A_\phi \in \cat{Bimod}(A)$, which is $A$ as a left $A$-module, but for which the right $A$-action is twisted by $\phi$. 
But $\cat{Bimod}(A) \simeq \Vect$, so
 there must be a bimodule isomorphism $A_\phi \isom A$. This map must take the unit $1 \in A_\phi$ to some (necessarily invertible) element $f\in A$, and compatibility with the left $A$-action says that it takes $a \mapsto af$. But compatibility with the right $A$-action then says that $a \phi(b) f= afb$ for all $a$, and so $\phi$ is inner. All together, the sequence
$$ \{1\} \to \bK^\times \to A^\times \to \Aut(A) \to \{1\}$$
is seen to be exact, and a central extension. The action $G \to \Aut(A)$ pulls back this central extension to an extension of $G$, classified by the anomaly $\alpha \in \H^2(G; \bK^\times)$.

In terms of the physical description from \S\ref{subsec:azumaya-physics}, the anomaly $\alpha$ arises as follows. The action of $G$ on $A$ determines an action of $G$ on the $(1+1)$-dimensional phase of matter built from $A$. Stack this $G$-equivariant phase of matter with a copy of the phase of matter built from $A^\op$ with trivial $G$-action. (``Stacking'' means to consider the noninteracting tensor product of the two physical systems. Physically, you place them on separate sheets separated by an infinitely strong insulating material, and then press the combined solar panel into a single layer.) Ignoring the $G$-symmetry, we find the trivial phase of matter, as it is built from the Morita-trivial algebra $A \otimes A^\op$. But the $G$-equivariance is nontrivial: 
we have built a ``symmetry-protected trivial'' phase. Bosonic SPT phases in $1+1$ dimensions are classified by ordinary group cohomology. 
(In the condensed matter sense, that claim is a widely believed assertion rather than a theorem; c.f.\ \cite[\S5.1]{GJFa}. It is a theorem in the TQFT sense: it is equivalent to the abstract algebraic construction given above.)

There is a final construction of the anomaly $\alpha$ which is useful to describe, since it is the easiest to generalize to the VOA setting. (We will describe it assuming the action of $G$ on $A$ is faithful; it is not hard to extend the construction to the non-faithful case.) Consider the subalgebra $A^G \subset A$ fixed by the $G$-action. Provided $G$ is finite and $\bK$ has characteristic $0$, this subalgebra remains semisimple. Its \define{commutant} is:
$$ B := \{b \in A \st [b,a] = 0 \;  \forall a \in A^G\}.$$
We claim that $B$ is a twisted group algebra for $G$. To prove this, for each $g \in G$ set:
$$ B_g := \{b \in A \st ba = g(a)b \; \forall a \in A\}.$$
Then $B_g$ is one-dimensional by the Skolem--Noether theorem (it is at most one-dimensional since $\cZ(A) \cong \bK$, and it is at least one-dimensional since we have already shown that all automorphisms of $A$ are inner) and $B_g \cap B_{g'} = \{0\}$ if $g \neq g'$ since the action is faithful. It is clear that $B_g B_{g'} = B_{gg'}$, and so $\bigoplus_{g\in G} B_g \subset A$ is a subalgebra. But the commutant of $\bigoplus_{g\in G} B_g$ in $A$ is  $A^G$, and so
$$ B = \bigoplus_{g\in G} B_g$$
by the double commutant theorem.

Finally, being a twisted group algebra, $B$ determines a cohomology class $\alpha \in \H^2(G; \bK^\times)$, easily identified with the cohomology class that classified the central extension discussed above.

\subsection{Gauging}\label{assoc-gauging}

Gauge anomalies arise as the obstruction to gauging symmetries. Consider first the case $A = \operatorname{Mat}_n(\bC)$, which corresponds via the bulk-boundary construction of \S\ref{subsec:azumaya-physics} to a quantum mechanics model with Hilbert space~$\bC^n$. A gaugeable group of symmetries is a group~$G$ acting on~$\bC^n$; then the gauged theory has as its Hilbert space the fixed subspace $(\bC^n)^G$. The anomaly must be trivialized in order to promote a group of automorphisms of $A$ to a group of automorphisms of the Hilbert space.

For associative algebras, gauging is also called ``quantum Hamiltonian reduction'';
the choice of trivialization of the anomaly corresponds to the ``quantum comoment map''.
Quantum Hamiltonian reduction can be written in various ways. Here is one.
Suppose that $G$ acts on $A$ with trivialized anomaly, meaning that we have chosen a way to lift  the homomorphism $G \to \Aut(A)$ to a homomorphism $G \to A^\times$. Then in particular the fixed subalgebra $A^G$ and the untwisted group algebras $\bK[G]$ are commutants in the Azumaya algebra $A$. This provides an isomorphism between their centres: $\cZ = \cZ(A^G) \cong \cZ(\bK[G])$, the algebra of class functions. The homomorphism $G \to \{1\}$ determines a homomorphism $\cZ(\bK[G]) \to \bK$. Specifically, $\cZ(\bK[G])$ arises as the algebra of $\bK$-valued global functions on the Deligne--Mumford stack $G/G$ corresponding to the adjoint action, and $\bK = \cO(\mathrm{pt}/G)$, and the map $\cZ(\bK[G]) \to \bK$ is restriction along the inclusion $\mathrm{pt}/G \to G/G$.
 The \define{gauging} of $A$ by $G$ is by definition the associative algebra:
\begin{equation}\label{assoc-gauging-eqn}
 A\sslash^\beta G := A^G \otimes_{\cZ} \bK.
\end{equation}
The superscript emphasizes that the gauging depends on the choice of trivialization $\beta$ of the anomaly $\alpha$.
The set of trivializations, if nonempty, is a torsor for $\H^1(G; \bK^\times)$.

It can happen that $A \sslash^\beta G$ is the zero algebra. If it is not zero, then it is Azumaya.
One way to see this is to realize $A \sslash^\beta G$ as the algebra of endomorphisms of the left $A$-module $A \otimes_G \bK$ (where $G$ acts on $A$ by right multiplication via the homomorphism $\beta : G \to A^\times$).
Indeed, such a description makes clear that if $A\sslash^\beta G \neq 0$, then $A$ and $A\sslash^\beta G$ are Morita equivalent (compare \cite{MR3581593}).

\subsection{Galois action and the cohomological Brauer group} \label{assoc-cohom-Brauer}

Given a $\bK$-algebra $A$ and a field automorphism $\gamma$ of $\bK$, there is a new $\bK$-algebra $^\gamma A$ which is equal to $A$ as a ring, but with $\bK$-structure twisted by $\gamma$. If $A$ is Azumaya then so is $^\gamma A$, and if $A$ has an action of $G$ by algebra automorphisms, then so does $^\gamma A$. There are no surprises about how anomalies transform under twisting by $\gamma$: if $G$ acts on $A$ with gauge anomaly $\alpha \in \H^2(G; \bK^\times)$, then it acts on $^\gamma A$ with anomaly~$\gamma(\alpha)$. In other words, the anomaly $\alpha$ has coefficients not just in $\bK^\times$-as-an-abelian-group, but in $\bK^\times$-as-a-Galois-module.

This Galois equivariance leads to a cohomological characterization of Azumaya algebras. Suppose that $A$ and $B$ are Azumaya algebras over $\bK$ and that $\bK \subset \bL$ is a Galois extension with Galois group $\Gal(\bL/\bK)$ such that $A \otimes_\bK \bL \cong B \otimes_\bK\bL$ as $\bL$-algebras; let $f : A \otimes_\bK \bL \to B \otimes_\bK \bL$ denote a choice of isomorphism. Then $A \otimes_\bK \bL$ carries two different $\Gal(\bL/\bK)$-actions, providing the different descent data: letting $\rho$ denote the standard action of $\Gal(\bL/\bK)$ on $\bK$, the two actions are $\rho_A: \gamma \mapsto \id_A \otimes \rho(\gamma)$ and $\rho_B : \gamma \mapsto f^{-1} \circ (\id_B \otimes \rho(\gamma)) \circ f$. Neither $\rho_A(\gamma)$ nor $\rho_B(\gamma)$ is $\bL$-linear, but the composition
$$ \tilde\kappa ( \gamma) := \rho_B(\gamma) \circ \rho_A(\gamma)^{-1} $$
is. This $\tilde\kappa$ is not a homomorphism from $\Gal(\bL/\bK)$ to $\Aut(A\otimes_\bK \bL)$, but it is a twisted cocycle --- $\rho_A$ determines an action of $\Gal(\bL/\bK)$ on $\Aut(A\otimes_\bK \bL)$ --- and so determines a twisted cohomology class 
$$ \kappa \in \H^1(\Gal(\bL/\bK); \Aut(A\otimes_\bK \bL)).$$
It is a basic theorem that this class, together with $A$, determines the isomorphism type of $B$. 

The earlier discussion of gauge anomalies amounts to the fact that there is a universal class $\alpha^{\text{univ}} \in \H^2(\Aut(A \otimes_\bK \bL); \bL^\times)$ classifying the extension $1 \to \bL^\times \to (A \otimes_\bK \bL)^\times \to \Aut(A\otimes_\bK\bL) \to 1$. Actually, our discussion did not quite imply this, because we worked only with finite groups $G$. But the image of the cocycle $\tilde\kappa$, together with its Galois-conjugates, is a finite subset of $\Aut(A \otimes_\bK \bL)$ which is easily seen to be closed under composition, and so a subgroup, and so the earlier discussion does imply that $\alpha^{\text{univ}}$, even if it did not exist on all of $\Aut(A \otimes_\bK \bL)$, does exist after being restricted to the image of $\kappa$.

The composition of $\alpha^{\text{univ}} \circ \kappa$ is a class in the twisted cohomology $\H^2(\Gal(\bL/\bK); \bL^\times)$. Given the Morita class of $A$, this twisted cohomology class determines the Morita equivalence class of $B$.

Finally, it is a fundamental fact that for every Azumaya algebra $B$ over $\bK$, there is a Galois extension $\bK \subset \bL$ such that $B \otimes_\bK \bL \cong \operatorname{Mat}_n(\bL)$ for some $n$. 
Fixing $A = \operatorname{Mat}_n(\bK)$ and taking the union over extensions $\bL \supset \bK$, one discovers:
$$ \Br(\bK) = \{\text{Azumaya algebras over $\bK$ up to Morita equivalence}\} \cong \H^2(\Gal(\bar{\bK}/\bK); (\bar{\bK})^\times).$$
The left-hand side is called the \define{Brauer group} of $\bK$ and the right-hand side is called the \define{cohomological Brauer group}. (In positive characteristic, $\bar\bK$ means the separable, rather than algebraic, closure.)

\section{Gauge anomalies for VOAs}\label{sec:VOAs}

In the remainder of this paper, we will repeat the story from \S\ref{sec:warmup} with associative algebras (the observables for $(0+1)$-dimensional quantum field theories) replaced by vertex operator algebras ($(1+1)$-dimensional quantum field theories). The question of how to gauge an action of a finite group on a VOA was first studied in \cite{MR1003430,MR1128130}. The mathematical construction of gauge anomalies is due to \cite{MR1923177}; this section consists primarily of an interpretation of that paper.

\subsection{Holomorphic VOAs} \label{voas-holo}

We will not review the definition of vertex operator algebra --- 
an excellent reference is \cite{MR2023933}
--- preferring instead to indicate the meanings of words used in VOA theory by analogizing with the associative algebra case. We first point out that the notion of VOA is purely algebraic: it makes sense over any field. 
We will work only over $\bC$ and its sub-fields, which  simplifies the story. (In positive characteristic, one must take care with divided powers. The  introduction of vertex algebras in \cite{MR843307} already included the positive-characteristic case.) In particular, if $V$ is a VOA over $\bK$, then so are its Galois conjugates $^\gamma V$ for $\gamma \in \Gal(\bK/\bQ)$. 
In this paper VOAs are always 
simple,
self-dual (i.e.\ admitting an invariant nondegenerate bilinear form),
and
$\bN$-graded by the action of $L_0$.

Given a VOA $V$, there are various categories of ``$V$-modules,'' with different regularity conditions placed on the module. We will focus on ``admissible'' modules (we will not give the definition), writing $\Rep(V)$ for the category thereof.
VOA modules correspond not to left-modules for an associative algebra $A$ but rather are analogous to $A$-$A$-bimodules: $\Rep(V)$ is the VOA analog of $\cat{Bimod}(A)$, and the Zhu algebra of $V$ (defined so that
the irreducible $V$-modules are in natural bijection with irreducible modules over the Zhu algebra; Zhu's construction in \cite{MR1317233} is technical, but see \cite{MR2568399} for a readable survey)
  is analogous to the enveloping algebra $A \otimes A^\op$.

The analog of semisimplicity for an algebra is called ``regularity'': a VOA $V$ is \define{regular} if it is $C_2$-cofinite and the category $\Rep(V)$ of admissible modules is semisimple with finitely many simples (a condition called ``rationality,'' but 
 this paper avoids that term so as not to conflict with the question of definability over $\bQ$).
Suppose that we are working over $\bK = \bC$. Then regularity implies that $\Rep(V)$ is a modular tensor category
\cite{MR2468370}.
The braided monoidal structure is the VOA analog of the fact that $\cat{Bimod}(A)$ is naturally a monoidal category.

The word analogous to ``Azumaya'' is ``holomorphic'': $V$ is \define{holomorphic} if $\Rep(V) \simeq \Vect$.

\subsection{Physical interpretation} \label{subsec:voa-physics}

It is believed that $(2+1)$-dimensional topological phases of matter are classified by MTCs together with some mild extra data about the central charge. (Namely, each MTC has a central charge $c$ living in $\bQ/8\bZ$; the extra data is a lift to $\bQ$.  The root of unity $\alpha = \exp(2 \pi i c/8)$ is sometimes called the \define{(gravitational) anomaly} of the MTC.) There is no proof of this belief, largely due to the lack of a definition of ``$(2+1)$-dimensional topological phase of matter,'' but a TQFT version is proven in \cite{BDSPV}, where it is confirmed that the Reshetikhin--Turaev construction \cite{MR1091619} gives a precise classification of once-extended 3d TQFTs. If $V$ is a regular VOA over $\bC$, then the MTC $\Rep(V)$ determines a $(2+1)$-dimensional topological phase of matter, and the VOA $V$  arises as the algebra of observables for a distinguished $(1+1)$-dimensional chiral conformal boundary condition.

If $V$ is holomorphic, then the bulk $(2+1)$-dimensional phase of matter is almost trivial: it is a short range entangled (SRE) phase, encoding only the central charge of $V$. 
(The central charge corresponds to the ``beyond group cohomology layer'' in \cite[\S5.2]{GJFa}.) 
Consider a bulk-boundary system with $V$ determining the boundary physics and with bulk physics this SRE phase. Since the bulk phase is SRE, any choice of local boundary coordinates identifies the boundary physics with a true quantum field theory, but the theory isn't truly conformal: conformal transformations change the partition function of the theory in a simple way determined by $c$, and so $c$, or equivalently the bulks SRE phase, encodes the ``gravitational anomaly'' of the holomorphic CFT described by $V$.

\subsection{Construction of gauge anomalies} \label{voa-gauge-construction}

Let $V$ be holomorphic VOA defined over $\bC$ and let $G \subset \Aut(V)$ be a finite group of VOA automorphisms of $V$. We will associate to these data a \define{gauge anomaly} $\alpha \in \H^3(G;\bC^\times)$. We will do so by imitating the construction from \S\ref{assoc-construction}. There is no VOA analog of the exact sequence (\ref{eqn-assoc}). But we can consider the analog of the commutant subalgebras $A^G \subset A \supset B \cong \bK^\alpha[G]$.

First, it makes perfect sense in VOA theory to consider the $G$-fixed subalgebra $V^G \subset V$. 
There is a version of commutant subalgebras that can be said entirely within the language of vertex algebras, called the ``coset construction'' \cite{MR778819,MR1159433}: the \define{coset} of a sub-vertex-algebra $W \subset V$ is the sub-vertex-algebra $W' \subset V$ consisting of all elements of $V$ that have trivial operator product expansion with all elements of $W$. 
In our case, however, $V^G \subset V$ is a conformal embedding (i.e.\ the image of $V^G$ contains the Virasoro vector of $V$) 
and so, since $V$ is assumed simple,  $(V^G)'$ is the trivial VOA.

In order to have a meaningful commutant of $V^G \subset V$, we must leave the world of VOAs, which encode the vertex operators in a chiral CFT, and also consider topological line operators. One expects on physical grounds that topological line operators for a fixed CFT form a category, because there can be topological junctions between line operators, and that they can be fused, providing the category of line operators with a monoidal structure. (A typical CFT will have infinitely many nonisomorphic simple topological line operators and so the category of line operators is not fusion. C.f.\ \cite{Andre-bicommutant}, which uses the term ``soliton'' for what we call ``topological line operator.'')

For each $g\in G$, consider the category $\cC_g$ of topological line operators $X$
 such that, for each vertex operator $v \in V$, 
 $$ 
 \begin{tikzpicture}[baseline=(v.base)]
   \draw[ultra thick] (0,0) -- (0,2); \path (.25,.25) node {$X$}; \path (1,1) node (v) {$\bullet \,v$};
   \draw[ultra thin] (-2,0) -- (2,0) -- (2,2) -- (-2,2) -- cycle;
 \end{tikzpicture}
 \quad=\quad
 \begin{tikzpicture}[baseline=(v.base)]
   \draw[ultra thick] (0,0) -- (0,2); \path (.25,.25) node {$X$}; \path (-1,1) node (v) {$\bullet \,g(v)$};
   \draw[ultra thin] (-2,0) -- (2,0) -- (2,2) -- (-2,2) -- cycle;
 \end{tikzpicture}
  $$
 in the sense that if you take the operator-valued holomorphic function ``insert the vertex operator $v$ at a point to the right of the line operator $X$'' and analytically continue it over the $X$-line, you get the operator-valued holomorphic function ``insert the vertex operator $g(v)$ at a point to the left of the line operator $X$.''

There is a natural functor from this category $\cC_g$ to what in the VOA literature is called the category of ``$g$-twisted $V$-modules.''
Indeed, a typical topological line operator $X$ in a CFT does not admit any topological ``endings,'' i.e.\ there are usually no topological junctions between $X$ and the trivial line operator. But there is a vector space of chiral conformal endings of $X$. (Compare: The topological endings of the trivial line operator are the vertex operators in $V$ of conformal dimension $0$, which is to say just the multiples of the vacuum. The chiral conformal endings of the trivial line are all of $V$.) When $X$ is a line operator satisfying $Xv = g(v)X$, the space of chiral conformal endings for $X$ is a $g$-twisted module.

The functor $\cC_g \to \{\text{$g$-twisted $V$-modules}\}$, $X \mapsto \{\text{endings of $X$}\}$, is expected to be an equivalence, and so we will conflate the two categories. The only way it could fail to be an equivalence is if there are line operators in $\cC_g$ which admit no endings at all. Such line operators can exist in a TQFT (assuming the phrase ``TQFT'' implicitly restricts just to topological endings of topological line operators) but should not exist for a CFT described by a VOA.

When $V$ is holomorphic, the category $\cC_g$ is (noncanonically) equivalent to $\Vect$ \cite[Theorem 2]{MR1794264}.
Consider the direct sum of categories $\cC = \bigoplus \cC_g$. Under the physical expectation that line operators can be fused, $\cC$ will be a $G$-pointed fusion category.
The nontrivial data of a $G$-pointed fusion category is its associator; up to equivalence of $G$-pointed fusion categories, that data is a cohomology class $\alpha \in \H^3(G;\bC^\times)$ \cite[Example 2.3.8 and Proposition 2.6.1]{EGNO}. By construction $V^G$ is the commutant VOA of $\cC \subset V$, since it consists of those vertex operators that commute with the line operators in $\cC$,
and one expects a ``double commutant theorem'' realizing $\cC = \bigoplus \cC_g$ as the commutant of $V^G$.

Actually, working purely in VOAs, it is not known that one can fuse arbitrary twisted modules. We will therefore impose the following widely-believed assumption:

\mytheorem[Regularity Assumption]{Let $V$ be a holomorphic VOA and $G \subset \Aut(V)$ a finite group. Then the fixed sub-VOA $V^G$ is regular.}

\noindent The main result of \cite{CarnahanMiyamoto} implies that this Regularity Assumption holds if $G$ is solvable.
With the Regularity Assumption in place, the results of \cite{MR1923177} fully justify the story told in this section: in particular, the category $\cC$ is fusion, and is the commutant of $V^G$, and so the gauge anomaly $\alpha \in \H^3(G;\bC^\times)$ is defined.

\subsection{Gauging}\label{voa-gauging}

To justify the name ``gauge anomaly,'' we must construct, for each trivialization~$\beta$ of~$\alpha$, a \define{gauging} $V\sslash^\beta G$ of $V$. In the physicist's style we will first give a recipe for $V\sslash^\beta G$. In \S\ref{voa-linalg} we will explain the question that this recipe answers.
The recipe will imitate the quantum Hamiltonian reduction formula~(\ref{assoc-gauging-eqn})
from \S\ref{assoc-gauging}: for Azumaya algebras, the gauging was $A \sslash^\beta G := A^G \otimes_{\cZ} \cO(\mathrm{pt}/G)$, where $\cZ \cong \cZ(A^G) \cong \cZ(B)$ and where $\beta$ identified the twisted group algebra $B = \bigoplus B_g \cong \bK^\alpha[G]$ with the untwisted group algebra $\bK[G]$.

In the VOA case, the gauge anomaly $\alpha$ arises as the associator on the $G$-pointed fusion category $\cC = \bigoplus_g \cC_g$. The trivialization is therefore the data of an equivalence $\beta : \cC \simeq \Vect[G]$, where $\Vect[G]$ means the $G$-pointed fusion category of $G$-graded vector spaces with trivial associator. This is directly analogous with the story about $B \cong \bK^\alpha[G] 
\overset \beta
\cong \bK[G]$.

We need an analog for the centre $\cZ(A^G)$. The \define{centre} of any algebra is its commutant as a subalgebra of itself. For $W$ a regular vertex algebra, we could take the coset of $W \subset W$, but we would get a trivial answer. Instead, as in \S\ref{voa-gauge-construction}, we will take its categorical commutant, i.e.\ the category of the line operators on $W$ that commute with all vertex operators in $W$. The category of such line operators is precisely $\Rep(W)$: each line operator has a $W$-module worth of chiral conformal endings (and we take for granted the assumption that every line operator has endings, so that the functor from line operators to modules is an equivalence).

A second justification for setting $\cZ(W) := \Rep(W)$ comes from comparing \S\ref{subsec:azumaya-physics} and \S\ref{subsec:voa-physics}. The centre $\cZ(B)$ of a semisimple algebra $B$ arises as the algebra of operators for the $(1+1)$-dimensional topological order constructed from $B$, whereas $\Rep(W)$ are the operators for the $(2+1)$-dimensional topological order constructed from the regular VOA $W$.

The isomorphism $\cZ(A^G) \cong \cZ(B)$ used in \S\ref{assoc-gauging} has an analog: conditional on the Regularity Assumption, it is shown in \cite{MR1923177} that there is a canonical equivalence $\cZ(V^G) := \Rep(V^G) \simeq \cZ(\cC)$, where the right-hand side denotes the Drinfeld centre of the fusion category $\cC$.  

Moreover, analogous to the map $\cZ(\bK[G]) = \cO(G/G) \to \bK = \cO(\mathrm{pt}/G)$, there is a canonical map $\cZ(\Vect[G]) \to \Rep(G)$, coming from realizing $\cZ(\Vect[G])$ and $\Rep(G)$ as categories of sheaves on $G/G$ and $\mathrm{pt}/G$ respectively. The map $\cZ(\Vect[G]) \to \Rep(G)$ extends to an equivalence $\cZ(\Vect[G]) \simeq \cZ(\Rep(G))$ of MTCs \cite[Example 3.15]{MR3077244}.
One can recover the trivialization $\beta$ from the composition $\cZ := \cZ(\cC) \overset\beta\simeq \cZ(\Vect[G]) \to \Rep(G)$.

All together, we can analogize  equation~(\ref{assoc-gauging-eqn}) to
\begin{equation} \label{voa-gauging-eqn}
  V\sslash^\beta G := V^G \otimes_\cZ \Rep(G).
\end{equation}
The last step is to explain how to tensor a VOA with a fusion category. The short answer is that $V\sslash^\beta G$ is a conformal extension of $V^G$ determined by $\cZ \to \Rep(G)$: the map $V^G \mono V\sslash^\beta G$ is analogous to the map $A^G \epi A \sslash^\beta G := A^G \otimes_\cZ \bK$; in categories, there is no real distinction between quotients and extensions (e.g.\ extensions of algebras often lead to quotients of module categories).

The long answer is that the data of the map $\cZ \to \Rep(G)$ is equivalent to the data of an ``\'etale algebra object'' $E \in \cZ$, i.e.\ a separable braided-commutative algebra which contains only one copy of the unit object of $\cZ$. The equivalence identifies $\Rep(G)$ as the category of $E$-module-objects in $\cZ$, and identifies $E$ as the $\cZ$-object worth of ``internal endomorphisms'' of the unit in $\Rep(G)$.
But \'etale algebra objects in $\Rep(V^G)$ are precisely conformal extensions of $V^G$.
All together, equation~(\ref{voa-gauging-eqn}) says that $V\sslash^\beta G$ is the conformal extension of $V^G$ corresponding to the \'etale algebra $E \in \Rep(V^G)$ corresponding to the map $\Rep(V^G) \overset\beta\simeq \cZ(\Vect[G]) \to \Rep(G)$.

The \'etale algebra object $E \in \cZ$ is ``lagrangian'' in the sense that it realizes $\cZ$ as the Drinfeld centre of $\Rep(G)$. It follows that $V \sslash^\beta G$ is holomorphic.

\subsection{Linear algebraic description}\label{voa-linalg}

The previous sections provided a recipe that took in a holomorphic VOA $V$ and a finite group $G\subset \Aut(V)$ and a trivialization $\beta$ of the anomaly and produced (assuming the Regularity Assumption) a VOA $V \sslash^\beta G$ that we called the ``gauging'' of $V$ by $G$. We did not explain abstractly what ``gauging'' means, nor did we prove that one needs precisely a trivialization of $\beta$. We will do so now. Specifically, we will arrive at the answer by explaining how $V \sslash^\beta G$ decomposes as a module over its subalgebra $V^G$. We call this a ``linear algebraic description'' of $V \sslash^\beta G$ since we will not try to describe the operator product.

The first step is to describe the category $\Rep(V^G)$ in which $V \sslash^\beta G$ lives. Suppose that $G$ acts with anomaly $\alpha \in \H^3(G;\bC^\times)$, which we don't yet assume to be trivializable; then $\Rep(V^G) \simeq \cZ(\Vect^\alpha[G])$. 
Slightly abusively, let us write $G/G$ for the set of conjugacy classes in $G$, and for each $g \in G/G$ choose a representative also called $g$.  Let $C(g) \subset G$ denote its centralizer. There is a \define{slant product} map $\iota_g : \H^3(G;\bC^\times) \to \H^2(C(g);\bC^\times)$, given on cocycle representatives by 
$$ (\iota_g\alpha)(x,y) = \alpha(g,x,y) - \alpha(x,g,y) + \alpha(x,y,g). $$
(Here and throughout, we write the group law on  $\H^\bullet(G;\bC^\times)$ additively.)

Write $\Rep^{\iota_g\alpha}(C(g))$ for the $\iota_g\alpha$-block in the category of projective $C(g)$-modules.
Then, as a linear category,
$$ \cZ(\Vect^\alpha[G]) \simeq \bigoplus_{g\in G/G} \Rep^{\iota_g\alpha}(C(g)).$$
For each $g \in G/G$, let $V_g \in \cC_g$ denote (a choice for) the simple $g$-twisted $V$-module. Then $V_g$ is naturally a $\iota_g\alpha$-projective $C(g)$-module, i.e.\ $V_g \in  \Rep^{\iota_g\alpha}(C(g))$. 
In fact, under the equivalence $\Rep(V^G) \simeq \cZ(\Vect^\alpha[G])$, the simple objects of $\Rep^{\iota_g\alpha}(C(g))$ are identified with the $V^G$-submodules of $V_g$. 

Suppose now that we choose a trivialization $\beta$ of $\alpha$. Then $\iota_g\beta$ (represented by the cocycle $(\iota_g\beta)(x) = \beta(g,x) - \beta(x,g)$) is a trivialization of $\iota_g\alpha$, i.e.\ it determines a one-dimensional $\iota_g\alpha$-projective $C(g)$-module that we will call $\bC_{\iota_g\beta} \in \Rep^{\iota_g\alpha}(C(g))$. Given $M \in \Rep^{\iota_g\alpha}(C(g))$, let $M^{\iota_g\beta} := \hom(\bC_{\iota_g\beta}, M)$. It is the space of ``fixed points'' of the $C(g)$-action if we use $\iota_g\beta$ to identify $\Rep^{\iota_g\alpha}(C(g))$ with $\Rep(C(g))$.

Then, as a representation of $V^G$, we have
\begin{equation} \label{voa-linalg-eqn}
 V \sslash^\beta G \cong \bigoplus_{g \in G/G} V_g^{\iota_g\beta}. 
\end{equation}
One can show this by inspecting the lagrangian \'etale algebra object in $\Rep(\Vect[G])$ corresponding to the fusion category $\Rep(G)$.

An immediate consequence of equation~(\ref{voa-linalg-eqn}) is that the only common $V^G$-submodule of $V$ and $V \sslash^\beta G$ is $V^G = V_e^{\iota_e\alpha}$ itself (where of course $e\in G$ denotes the unit). Based on this, we take:

\mytheorem[Definition]{Let $V$ be a holomorphic VOA and $G \subset \Aut(V)$ a finite group of automorphisms such that the Regularity Assumption holds. A \define{gauging of $V$ by $G$} is any holomorphic VOA $V\sslash G$ extending the fixed subalgebra $V^G$ such that the only common $V^G$-submodule shared by $V$ and $V\sslash G$ is $V^G$ itself. The action of $G$ on $V$ is \define{nonanomalous} if a gauging exists.}

We have seen that if the anomaly admits a trivialization $\beta$, then $V \sslash^\beta G$ is a gauging in this sense.
Conversely, by studying lagrangian algebras in $\cZ(\Vect^\alpha[G])$, one can show that the action of $G$ on $V$ is nonanomalous if and only if the anomaly is trivializable, and that gaugings are in bijection with trivializations.
\vspace*{-3pt}

\section{Main result}\label{sec:mainthm}

Most of \S\ref{sec:VOAs} requires working over the complex numbers $\bC$ for a simple reason: there is no purely algebraic way to talk about ``line operators,'' because there is no purely algebraic way to talk about real submanifolds of $\bC$. Said another way, the MTC structure on $\Rep(W)$, for $W$ a regular VOA, requires {a priori} that we work over $\bC$: algebraically, one can contemplate intertwining operators parameterized by configurations of points on $\bP^1$; only after taking $\bC$-points can one identify configuration spaces on $\bP^1$ with classifying spaces for braid groups. Indeed, the construction of the MTC structure on $\Rep(W)$ given in \cite{MR2468370} repeatedly chooses parameters living on real intervals in $\bC$.

Thus the  definition of  the anomaly $\alpha$ of an action of $G$ on a holomorphic VOA $V$ requires that $V$ be defined over $\bC$. But suppose that $\gamma$ is a field automorphism of $\bC$. Then there is a Galois-conjugate VOA $^\gamma V$ equal to $V$ as a VOA over $\bQ$ but with the $\bC$-structure twisted by $\gamma$. It is also holomorphic, since its linear category of representations is defined algebraically, and also has a $G$-action, hence its own anomaly.
Our goal in this section is to prove:

\mytheorem[Main Theorem]{
Let $V$ be a holomorphic VOA over $\bC$ and $G \subset \Aut(V)$ a finite group of automorphisms, and let $\gamma$ be a field automorphism of $\bC$ and $^\gamma V$ the Galois-conjugate VOA. Conditional on the Regularity Assumption, let $\alpha \in \H^3(G;\bC^\times)$ denote the anomaly of the action of $G$ on $V$. Then the anomaly of the action of $G$ on $^\gamma V$ is $\gamma^2(\alpha)$.
}\vspace*{-7pt}

\subsection{Recollection on the cyclotomic Galois group} \label{galois-recollections}

Let $\bmu \subset \bC^\times$ denote the group of roots of unity. We will write $\bmu$ additively; it is isomorphic to $\bQ/\bZ$ via the exponential map $x \mapsto \exp(2\pi i x)$. If $G$ is a finite group, then the map $\H^\bullet(G;\bmu) \to \H^\bullet(G;\bC^\times)$ is an isomorphism for $\bullet \geq 1$.

If $\gamma \in \Aut(\bC)$ and $x \in \bmu$ has order $N$, then $\gamma(x) \in \bmu$ also has order $N$. It follows that $\gamma$ acts on $\bmu$ by multiplication by some $n \in \widehat{\bZ}^\times$, where $\widehat{\bZ}$ denotes the profinite completion of $\bZ$ and $\widehat{\bZ}^\times$ is its group of units.
The map $\gamma \mapsto n : \Aut(\bC) \to \widehat{\bZ}^\times$ is a surjection of $\Aut(\bC)$ onto the cyclotomic Galois group $\Gal(\bQ^{\mathrm{cyc}}/\bQ) \cong \Aut(\bmu) \cong \widehat{\bZ}^\times$.

\subsection{Nonanomalous actions} \label{nonanom}

There are two purely algebraic parts of \S\ref{sec:VOAs}. First, the definitions in \S\ref{voas-holo} of the words ``holomorphic'' and ``regular'' refer only to the structure of the linear category of representations, and so if $V$ is holomorphic, then $^\gamma V$ is as well. Second, the definition of a ``gauging'' of $V$ by $G$ given in \S\ref{voa-linalg} is purely algebraic. 
To see this, first note that $(^\gamma V)^G = {^\gamma (V^G)}$, and so the expression $^\gamma V^G$ is unambiguous.
Second, suppose that  $V\sslash G$ is a gauging of $V$ by $G$, i.e.\ a holomorphic extension of $V^G$ which, as a $V^G$-module, shares only the submodule $V^G$ with~$V$. Then $^\gamma(V\sslash G)$ is a gauging of $^\gamma V$ by $G$, since the only common $^\gamma V^G$-submodule shared by $^\gamma V$ and $^\gamma(V\sslash G)$ is $^\gamma V^G$. In summary:

\mytheorem[Lemma]{If the action of $G$ on $V$ is nonanomalous, then so is the action of $G$ on $^\gamma V$. \qed}

\pagebreak[0]
It is not hard to prove that if $G$ acts on $V_1$ with anomaly $\alpha_1$ and on $V_2$ with anomaly $\alpha_2$, then it acts on $V_1 \otimes V_2$ with anomaly $\alpha_1 + \alpha_2$, where we have written the group law on $\H^3(G;\bC^\times)$ additively. 

\mytheorem[Corollary]{Suppose that $G$ acts on $V_1$ with anomaly $\alpha$ and on $V_2$ with anomaly $-\alpha$, and suppose that the Main Theorem holds for $V = V_1$. Then it also holds for $V = V_2$. \qed}

\subsection{A construction of Evans and Gannon} \label{egconst}

Conditional on the Regularity Assumption, the paper \cite{EvansGannon} constructs, for each finite group $G$ and each class $\alpha \in \H^3(G;\bmu)$, a VOA $V$ with a faithful $G$-action with anomaly $\alpha$. Given the corollary in \S\ref{nonanom}, it suffices to prove the Main Theorem for the VOAs constructed by Evans and Gannon. This section reviews their construction.

Let $G$ be a finite group. Then there is a finite abelian extension $1 \to A \to E \to G \to 1$ such that the restriction map $\H^3(G;\bmu) \to \H^3(E;\bmu)$ is the zero map. This fact has been rediscovered a number of times, most recently in \cite{PhysRevX.8.031048}. 
The Lyndon--Hochschild--Serre spectral sequence enhances the usual 5-term inflation-restriction exact sequence (for $\bmu$-cohomology) to one of the form:
\begin{multline*}
  1 \to \H^1(G;\bmu) \to \H^1(E;\bmu) \to H^1(A;\bmu) \to \H^2(G;\bmu) \to \H^2(E; \bmu) \\ \to \Lambda(G;A) \overset\delta\to \H^3(G; \bmu) \overset0\to \H^3(E;\bmu) \to \dots
\end{multline*}
where $\Lambda(G;A)$ maps naturally to $\H^2(A;\bmu)$, with kernel $\H^1(G; \H^1(A;\bmu)$.
It follows that for each $\alpha \in \H^3(G;\bmu)$, we can choose $\beta \in \Lambda(G;A)$ with $\delta\beta = \alpha$. 

The strategy of \cite{EvansGannon}, then, is to start with a holomorphic VOA $W$ with a faithful action by $E$ with trivialized anomaly. (They prove that the permutation action on $W = (V^\natural)^{\otimes |E|}$, where $V^\natural$ denotes the moonshine VOA, works. One can always arrange a nonanomalous $E$-action by starting with some holomorphic VOA $W$ on which $E$ acts with anomaly of order $N$ and then taking $W' = W^{\otimes N}$. It could happen that $W'$ is nonanomalous but not canonically so: the anomaly may not have a distinguished trivialization. But all choices of trivialization of the anomaly for $W'$ will lead to the same trivialization of the anomaly for $W'' = (W')^{\otimes |\H^2(E;\bmu)|}$.) 
 In particular the induced $A$-action on $W$ is nonanomalous, and so gaugings of $W$ by $A$ are parameterized by classes in $\H^2(A;\bmu)$. The central result underlying the construction is:

\mytheorem[Proposition (\cite{EvansGannon})]{
Conditional on the Regularity Assumption,
$W \sslash^\beta A$ carries a faithful $G$-action with anomaly $\alpha = \delta\beta$. \qed
}

In the Proposition, the VOA $W \sslash^\beta A$ depends only on the image of $\beta \in \Lambda(G;A)$ under the map $\Lambda(G;A) \to \H^2(A;\bmu)$, but the information in $\beta$ in the kernel of this map is used to determine the action of $G$. (The group $\H^1(A;\bmu)$ acts on any orbifold $W \sslash A$, and so an action of $G$ on $W$ can be adjusted by an element of $\H^1(G; \H^1(A;\bmu)) = \hom(G, \H^1(A;\bmu))$; compare the ``electromagenetic duality'' from \cite{FT18} or the ``T-duality'' from \cite{MR3916985}.)

\subsection{Completion of the proof}

Thus to prove the Main Theorem, it suffices to study the case $V = W\sslash^\beta A$ from \S\ref{egconst}. Specifically, we know that $G$ acts on $W\sslash^\beta A$ with anomaly $\alpha = \delta \beta$. Suppose that $\gamma \in \Aut(\bC)$ acts on $\bmu$ by multiplication by $n \in \widehat{\bZ}^\times$. (As in the previous sections, we will write $\bmu$ additively.) Then we wish to show that $G$ acts on $^\gamma(W \sslash^\beta A)$ with anomaly $n^2\alpha$. 
Since $\delta : \Lambda(G;A) \to \H^3(G;\bmu)$ is linear,
it suffices to prove:

\mytheorem[Claim]{$^\gamma(W \sslash^\beta A) \cong {^\gamma W} \sslash^{n^2\beta} A$.}

To prove the claim, we begin with the linear algebraic description of $W \sslash^\beta A$ from \S\ref{voa-linalg}. Using that $A$ is an abelian group, equation~(\ref{voa-linalg-eqn}) says that as a $W^A$-module,
$$ W \sslash^\beta A \cong \bigoplus_{a \in A} W_a^{\iota_a\beta}. $$
and so
$$ ^\gamma(W \sslash^\beta A) \cong \bigoplus_{a\in A} {^\gamma(W_a^{\iota_a\beta})}.$$

We must now understand the Galois conjugate $W^A$-module $^\gamma(W_a^{\iota_a\beta})$. 
By assumption $\gamma$ acts by multiplication by $n \in \widehat{\bZ}^\times$ on $\bmu$, and so
$$ ^\gamma(W_a^{\iota_a\beta}) \cong (^\gamma(W_a))^{\iota_a n\beta}.$$

What is $^\gamma(W_a)$? It is not $(^\gamma W)_a$. Indeed, the definition of ``$a$-twisted $W$-module'' refers explicitly to the eigenvalues of $a$.
These eigenvalues (or rather their logarithms, so that we can work additively)
 get multiplied by $n \in \widehat{\bZ}^\times$ under Galois conjugation by $\gamma$. To compensate for this one must divide $a$ by $n$. The result is:
$$ ^\gamma(W_a) \cong (^\gamma W)_{n^{-1}a}.$$

Taking the direct sum over $a\in A$, we find:
\begin{equation}\label{eqn-aclaim}
 ^\gamma(W \sslash^\beta A) \cong \bigoplus_{a\in A} {^\gamma(W_a^{\iota_a\beta})} \cong \bigoplus_{a \in A} (^\gamma W)_{n^{-1}a}^{\iota_a n\beta} \cong \bigoplus_{a' \in A} (^\gamma W)_{a'}^{\iota_{na'} n\beta} \cong \bigoplus_{a' \in A} (^\gamma W)_{a'}^{n^2\iota_{a'} \beta}.\end{equation}
Here we have reindexed $a' = n^{-1}a$. We also used that $\iota_{nx} \beta = n\iota_x \beta$. To see this, note that, for fixed $y$, $\iota_x\beta(y) = \beta(x,y) - \beta(y,x) = -\iota_y\beta(x)$ is a 1-cocycle  in the $x$ variable.

Equation~(\ref{eqn-aclaim}) establishes the Claim at the linear algebraic level, i.e.\ as modules over $^\gamma W^A$. But, since $A$ is an abelian group, \'etale algebras in $\cZ(\Vect[A]) \simeq \Rep(^\gamma W^A)$ are determined by their underlying modules. Thus the linear algebraic isomorphism implies a VOA isomorphism. 
Finally, the VOA isomorphism is compatible with the $G$-action.
This completes the proof of the Claim and hence of the Main Theorem.

\section{Questions and conjectures}\label{sec:conclusions}

\subsection{Moonshine for every group} \label{moonshine-speculation}

Since $\Aut(\bmu) \cong \widehat{\bZ}^\times$ is abelian, its action on $\bmu$ can be ``twisted'' by any power. Given $i\in \bZ$, we will write $\bmu^{\otimes i}$ for the abelian group $\bmu$ made into a $\widehat{\bZ}^\times$-, and hence $\Aut(\bC)$-, module via the action $n \triangleright x := n^i x$.
The Main Theorem can then be summarized as saying that as an $\Aut(\bC)$-module, gauge anomalies of VOAs live in $\H^3(-;\bmu^{\otimes 2})$.

\mytheorem[Corollary of the Main Theorem]{
Suppose that the action of $G$ on $V$ preserves a $\bK$-form for some field $\bK \subset \bC$. Then the anomaly of the action lives in
$
  \H^0\bigl(\Gal(\bar{\bK}/\bK); \H^3\bigl(G; \bmu^{\otimes 2}\bigr)\bigr).
$
}

Consider the case $\bK=\bQ$. The ``defining property of $24$'' from \cite{MR554399} --- that $n^2 \equiv 1 \mod 24$ for $n \in \widehat{\bZ}^\times$ --- implies that for any group $G$, $\H^0(\Gal(\bar{\bQ}/\bQ); \H^\bullet(G; \bmu^{\otimes 2}))$ is entirely $24$-torsion. Let $V^\natural$ denote the moonshine VOA. The monster group $\bM$ acts on $V^\natural$ preserving a $\bQ$-form~\cite{MR2928458}. Taking for granted the Regularity Assumption, let $\omega^\natural$ denote the anomaly of this action. In unpublished work \cite{Mason-Monster}, Mason computed the image of $\omega^\natural$ under the total slant product map
$$ \prod_{g\in \bM} \iota_g : \H^3(\bM;\bmu) \to \prod_{g \in \bM} \H^2(C(g); \bmu)$$
and showed that the image has exact order $24$. Together with the above corollary, we find:

\mytheorem[Corollary]{The anomaly $\omega^\natural$ of the action of $\bM$ on $V^\natural$ has exact order $24$.}

The exact order of $\omega^\natural$ was first computed in \cite{MR3916985} using a different argument;
that paper modelled CFTs in terms of conformal nets (rather than VOAs), for which the Regularity Assumption is known \cite{MR1806798} to hold.

The construction from \cite{EvansGannon}, reviewed in \S\ref{egconst}, establishes that there is a weak form of ``moonshine'' for every finite group $G$ and every anomaly $\alpha \in \H^3(G;\bmu)$. The resulting VOA
 is by no means canonical, and has astronomically large central charge (of order $|G| \times |\H^3(G;\bmu)|$). Various authors have speculated versions of the following (I heard it from I.\ Frenkel), while knowing that, as stated, it is surely too naive:

\mytheorem[Conjecture]{
For each quasisimple group $G$, there is a distinguished ``moonshine'' VOA representation $V_G$ of $G$. For the Monster group, $V_\bM = V^\natural$, and for compact Lie groups, $V_G$ is the level-one WZW model. The VOAs $V_G$ 
can be used in the classification of finite simple groups analogously to the way Lie algebras are used in the classification of simple Lie groups.
}

The words ``level-one'' in the phrase ``level-one WZW model'' of a simple Lie group $G$ refers to the fact that the anomaly for these models is a generator of the infinite cyclic group $\H^4(BG;\bZ)$. Theorem 8.6 of \cite{MR3990846} comes close to showing that $\omega^\natural$ generates $\H^3(\bM;\bmu) \cong \H^4(B\bM;\bZ)$. For all but finitely many finite simple groups $G$ of Lie type, $\H^3(G;\bmu)$ is cyclic, and the calculations in \cite{MR3990846} suggest that cyclicity also holds for ``most'' of the sporadic simple groups. Thus the above Conjecture may be refined by speculating that one way $V_G$ should be distinguished is that its anomaly should generate $\H^3(G;\bmu)$.

The calculations in \cite{MR3990846} also show that for ``most'' sporadic simple groups $G$, $\H^3(G;\bmu)$ is entirely 24-torsion. (The known exceptions are the Janko groups $J_1$, $J_2$, and $J_3$.) Given our Main Theorem, this increases the odds that the distinguished representation $V_G$ posited by the Conjecture might be defined over $\bQ$. For comparison, according to unpublished work of J.~Grodal~\cite{Grodal}, if $G$ is a quasisimple group of Lie type defined over $\bF_q$, then $\H^3(G;\bmu) \approx \bZ/(q^2-1)$. (The ``$\approx$'' sign indicates that there are some systematic adjustments needed to accommodate
twisted forms and
Schur multipliers of $G$, and that there are finitely many groups $G$ for which further corrections are necessary.) Note that, letting  $\bmu \subset \bar{\bF}_q^\times$ now denote the group of roots of unity in positive characteristic,
 $$\bZ/(q^2-1) \cong \H^0(\Gal(\bar{\bF}_q/\bF_q); \bmu^{\otimes 2}).$$
This is consistent with the construction of \cite{MR3456026} of a VOA over $\bF_q$ for each group of Lie type defined over $\bF_q$.
Missing is a good theory of ``gauge anomalies'' in positive characteristic.

\subsection{A vertex Brauer group}

Recall the story of the cohomological Brauer group outlined in~\S\ref{assoc-cohom-Brauer}. When a group $G$ acts on an Azumaya algebra $A$ over $\bar{\bK}$, the anomaly lives not just in $\H^2(G;\bar{\bK}^\times)$ thought of as an abstract group: under the $\Gal(\bar{\bK}/\bK)$-action on $\{$Azumaya algebras over $\bar{\bK}\}$, the anomaly transforms according to the action on (the coefficients of) $\H^2(G;\bar{\bK}^\times)$. This control over the Galois action implies that there is an anomaly-type map for Galois-twisted actions, which in the case of descent data for $\Gal(\bar{\bK}/\bK)$ provides a map
$$ \{\text{Azumaya algebras over $\bK$}\} \to \H^2\bigl(\Gal(\bar{\bK}/\bK); \bar{\bK}^\times\bigr).$$
Moreover, this map provides an isomorphism of the right-hand side with the Brauer group  $\Br(\bK)$.

Our Main Theorem suggests that there may be a similar story for VOAs. For $G$ a finite group, our Main Theorem identifies anomalies of actions of $G$ on holomorphic VOAs as living in $\H^3(G;\bmu^{\otimes 2})$. This suggests:

\mytheorem[Conjecture]{
For each field $\bK$, there is a \define{vertex Brauer group} of Morita equivalence classes of holomorphic VOAs over $\bK$ modulo tensoring with $E_{8,1}$. 
When $\bK$ is of characteristic $0$, then the vertex Brauer group is isomorphic to the Galois cohomology group $\H^3(\Gal(\bar{\bK}/\bK); \bmu^{\otimes 2})$.
}
We restrict to charactersitic zero for two reasons. First, there is so far no theory of 't Hooft anomalies of actions on VOAs in positive characterstic. Second, in the Azumaya case, the map
$$ \H^2\bigl(\Gal(\bar{\bK}/\bK);\bmu\bigr) \to \H^2\bigl(\Gal(\bar{\bK}/\bK);\bar{\bK}^\times\bigr) \cong \Br(\bK)$$
is an isomorphism only in characteristic $0$; in characteristic $p >0$ it is an isomorphism onto the prime-to-$p$ part of the right-hand side (see e.g.\ \cite[Chapter VI, Corollary 6.3.6]{MR2392026}).

We now outline, in parallel with \S\ref{assoc-cohom-Brauer}, how one could hope to prove the Conjecture. 
Suppose that $V$ and $W$ are holomorphic VOAs over $\bK$ such that $V \otimes_\bK \bar{\bK} \cong W \otimes_\bK \bar{\bK}$. If $V$ is fixed, then the data of the $\bK$-form $W$ is equivalent to the data of a $1$-cocycle with twisted coefficients 
$$\kappa \in \H^1\bigl(\Gal(\bar{\bK}/\bK); \Aut_{\bar{\bK}}(V \otimes_\bK \bar{\bK})\bigr).$$
Suppose that the  't Hooft anomalies of finite group actions on $V \otimes_\bK \bar{\bK}$ are the restriction of a universal anomaly living on $\Aut_{\bar{\bK}}(V \otimes_\bK \bar{\bK})$; according to our Main Theorem, this universal anomaly should live in
$$ \alpha \in \H^3\bigl( \Aut_{\bar{\bK}}(V \otimes_\bK \bar{\bK}); \bmu^{\otimes 2}\bigr).$$
Then the composition $\alpha \circ \kappa$ would be a twisted cohomology class
$$ \alpha \circ \kappa \in \H^3\bigl( \Gal(\bar{\bK}/\bK); \bmu^{\otimes 2}\bigr)$$
that depends only on the $\bK$-form $W$ of $V$.

Unfortunately, we lack direct access to this universal class $\alpha$: with current understanding of VOAs, we can only define anomalies for actions of finite groups, but $\Aut_{\bar{\bK}}(V \otimes_\bK \bar{\bK})$ is typically an infinite (in fact, affine algebraic) group. One could salvage this as follows. Suppose that there is a finite-degree Galois extension $\bK \subset \bL$ such that $V \otimes_\bK \bL \cong W \otimes_\bK \bL$. Then, still fixing $V$, the data of the $\bK$-form $W$ is given by a twisted cohomology class $\kappa_\bL \in \H^1(\Gal(\bL/\bK); \Aut_\bL(V \otimes_\bK\bL))$. Fix a cocycle for $\kappa_\bL$ and let $G \subset \Aut_\bL(V \otimes_\bK\bL)$ denote the union of the image of $\kappa_\bL$ together with its Galois conjugates. Then $G$ is a finite group, and so does have a well-defined anomaly $\alpha_\bL \in \H^0(\Gal(\bar{\bK}/\bL);\H^3(G; \bmu^{\otimes 2}))$.
Our proof of the Main Theorem does not provide sufficient control over the Galois action on cocycles --- at too many steps we worked instead simply with cohomology classes --- but it is reasonable to expect that there is a valid composition
$$ \alpha_\bL \circ \kappa_\bL \in \H^3\bigl( \Gal(\bL/\bK); \H^0\bigl(\Gal(\bar{\bL}/\bL); \bmu^{\otimes 2}\bigr)\bigr) \subset \H^3\bigl( \Gal(\bar{\bK}/\bK); \bmu^{\otimes 2}\bigr),$$
providing the class $\alpha \circ \kappa$ above.

Thus we expect that $\bK$-forms of a fixed holomorphic VOA lead to classes in the putative ``cohomological vertex Brauer group'' $\H^3(\Gal(\bar{\bK}/\bK); \bmu^{\otimes 2})$. In the associative  case, the next step is to classify Azumaya algebras over $\bar{\bK}$: they are all matrix algebras. In the vertex case, it is hopeless to try to classify holomorphic VOAs up to isomorphism, since in particular the set of even unimodular lattices injects into the set of holomorphic VOAs. (In small central charge, there is an expected classification \cite{MR1213740}, but it is not elementary.)

Instead, we should recognize the Azumaya classification as being about Morita equivalence classes: over $\bar{\bK}$, all Azumaya algebras are Morita equivalent. We can hope that a similar statement holds for VOAs. 
There is no consensus definition of ``Morita equivalence'' of VOAs. We will take the following definition:

\mytheorem[Definition]{Suppose that $V_1$ and $V_2$ are regular VOAs which are both conformal extensions of the same regular VOA $W \subset V_1,V_2$ (in particular,  all three of $V_1$, $V_2$, and $W$ have the same central charge). Then $V_1$ and $V_2$ are \define{Morita equivalent} if the corresponding \'etale algebras in $\Rep(W)$ are $E_2$-Morita equivalent.}
Implicit in the definition is the fact that, given a braided monoidal category $\cC$, there is a ``Morita $3$-category'' $\cat{Alg}_2(\cC)$ whose objects are braided-commutative algebra objects in $\cC$; compare \cite[\S8]{JFS}. In the case of a braided fusion category $\cC$, an $E_2$-Morita equivalence of \'etale algebras $A$ and $B$ is the same as a $\cC$-balanced Morita equivalence of fusion categories $\cat{Mod}(A) \simeq \cat{Mod}(B)$. It follows in particular that $\Rep(V_1) \simeq \Rep(V_2)$ as MTCs.

Suppose that $V_1$ and $V_2$ are holomorphic VOAs. Then, using that two fusion categories are Morita equivalent if and only if they have equivalent centres \cite[Theorem 3.1]{MR2735754}, one finds that $V_1$ and $V_2$ are Morita equivalent as soon as they share a common regular subalgebra.

The analog of the fact that all Azumaya algebras over $\bC$ are Morita trivial would then be the conjecture that the central charge is the only Morita invariant of a holomorphic VOA. Equivalently:

\mytheorem[Conjecture]{
Suppose $V_1$ and $V_2$ are two holomorphic VOAs over $\bC$ of the same central charge. Then there is a regular VOA $W$ of the same central charge embedding into both $V_1$ and $V_2$.
}
It is not at all clear whether to believe this conjecture, because we don't have enough understanding of the set of all VOAs of large central charge: the Conjecture holds for the few holomorphic VOAs that are known simply because the only  available methods to construct holomorphic VOAs are all of the form ``pass to a regular subalgebra, then extend to a holomorphic algebra in a different way'' and so always produce Morita-equivalent VOAs. 
(For instance, the Conjecture holds for all holomorphic VOAs of central charge $c\leq24$ by \cite{vEMS,MS19}. Holomorphic VOAs of larger central charge have not been well studied; the Conjecture applies, by construction, to those from \cite{largeC}.)
For lattice VOAs, the Conjecture follows from a $\bZ$-version of the Graham--Schmidt process, which produces, for any two unimodular lattices of the same rank, a common finite-index sublattice, and which provides the foundation of Kneser's neighbourhood method.
Physically, the Conjecture would follow from the expectation that every full conformal field theory with $c_L = c_R$ admits a gapped boundary condition.

Finally, if this Conjecture holds, then we could take the VOAs $E_{8,1}^{\otimes n}$ as our standard Morita-class representatives, where ``$E_{8,1}$'' means the level-one WZW theory for $E_8$ equipped with its ``Chevalley'' $\bQ$-form constructed in \cite{MR2928458}. This is why ``modulo tensoring with $E_{8,1}$'' appears in the statement of the Conjecture at the beginning of this subsection.

\subsection{K-theoretic interpretation of anomalies}

The existence of a vertex Brauer group isomorphic to $\H^3(\Gal(\bar{\bK}/\bK);\bmu^{\otimes 2})$ would imply the existence of a ``refined anomaly'' whenever the action of $G$ on a holomorphic vertex algebra $V$ preserves a $\bK$-form of $V$. Specifically, it would provide an anomaly map
\begin{equation}\label{eqn-veryrefinedanomaly}
 \{\text{holomorphic VOAs over $\bK$ with $G$-action}\} \to \H^3\bigl(G \times \Gal(\bar{\bK}/\bK); \bmu^{\otimes 2}\bigr).
\end{equation}
The ordinary gauge anomaly and the vertex Brauer class would be the restrictions to $\H^3(G;\bmu^{\otimes2})$ and $\H^3(\Gal(\bar{\bK}/\bK); \bmu^{\otimes 2})$ respectively.

The group $\H^3(G \times \Gal(\bar{\bK}/\bK); \bmu^{\otimes 2})$ is built, via a spectral sequence, from components of the form $\H^{3-i}(G; \H^{i}(\Gal(\bar{\bK}/\bK;\bmu^{\otimes 2})))$ for $0\leq i \leq 3$. Suppose that $\bK$ is a number field (i.e.\ a finite-degree extension of $\bQ$). Then the proposed vertex Brauer group is not itself particularly interesting: it follows from a theorem of Tate and Poitou \cite{MR0175892,PoitouLille} (see \cite[\S II.6.3, Theorem B]{MR1466966}) that $\H^3(\Gal(\bar{\bK}/\bK); \bmu^{\otimes 2})$ is isomorphic to $(\bZ/2)^N$, where $N$ is the number of field embeddings $\bK \mono \bR$. But the Galois cohomology groups $\H^i(\Gal(\bar{\bK}/\bK);\bmu^{\otimes 2})$ for $i\leq 2$ are quite interesting.
A theorem of Tate in the number-field case \cite{MR0422212,MR0429837}, and of Merkurjev and Suslin in general \cite{MR675529}, provides
  a close relationship between $\H^2(\Gal(\bar{\bK}/\bK);\bmu^{\otimes 2})$ and the second K-group $\rK_2(\bK)$:
\begin{equation}\label{MSthm}
 \H^2(\Gal(\bar{\bK}/\bK); \bmu^{\otimes 2}_n) \cong \rK_2(\bK) \otimes \bZ/n\bZ,\end{equation}
where $\bmu_n^{\otimes 2}$ denotes the $n$-torsion subgroup of $\bmu^{\otimes 2}$. 
(See \cite[Chapter 8]{MR2266528} and \cite[\S VI.4]{MR2392026}.)

Suppose that $G$ is a finite group acting on a holomorphic VOA $V$ over $\bK$.
 The refined anomaly~(\ref{eqn-veryrefinedanomaly}) would map to a class in $\H^1(G; \H^2(\Gal(\bar{\bK}/\bK);\bmu^{\otimes 2}))$. The isomorphism (\ref{MSthm}) then suggests:

\mytheorem[Conjecture]{
Each action of a finite group $G$  on a holomorphic VOA $V$ over $\bK$ determines 
 a homomorphism $G \to \rK_2(\bK)$.
} 

  Let $n$ denote the {exponent} of $G$, i.e.\ $g^n = 1$ for all $g \in G$. Then a map $G \to \rK_2(\bK)$ will land in the $n$-torsion subgroup of $\rK_2(\bK)$. Suppose that $\bK$ contains a primitive $n$th root of unity~$\zeta$. The choice of $\zeta$ determines an isomorphism between the $n$-torsion subgroup of $\rK_2(\bK)$ and the $n$-torsion subgroup of  $\mathrm{Br}(\bK)$, because it identifies both with the untwisted cohomology group $\H^2(\Gal(\bar{\bK}/\bK); \bZ/n)$. The choice of $\zeta$ is also what is needed to define algebraically the category $\cC_g$ of $g$-twisted $V$-modules.  Over an algebraically closed field, $\cC_g$ is (noncanonically) equivalent to $\Vect$, but such an equivalence need not survive Galois descent: rather, $\cC_g$ will be equivalent to the category of modules of some Azumaya algebra $A_g$.
  The map $g \mapsto [A_g] \in \Br(\bK)$ 
    presumably agrees with the map $G \to \rK_2(\bK)$ predicted in the Conjecture.

\subsection{Galois and Grothendieck--Teichmuller groups and modular data}

Let $W$ be a regular VOA over $\bC$. We emphasized at the start of \S\ref{sec:mainthm} that the MTC structure on $\Rep(W)$ depends, {a priori}, on the $\bC$-structure of $W$.
To control this {a priori} $\bC$-dependence one should study the following question. Choose $\gamma \in \Aut(\bC)$. As a $\bC$-linear category, $\Rep(^\gamma W) \simeq {^\gamma \Rep(W)}$, since the notion of ``representation'' is defined algebraic. How do the MTC structures compare? 

\mytheorem[Conjecture]{Fix a semisimple $\bC$-linear category $\cC$ with finitely many simples objects. The set of MTC structures on $\cC$ is permuted by a canonical action of $\Aut(\bC)$.}

Up to issues of completion, this action of $\Aut(\bC)$  should come from the injection of $\Gal(\bar{\bQ}/\bQ)$ into the Grothendieck--Teichmuller group, which is, again up to issues of completion, the automorphism group of the $E_2$ operad. MTCs are a particular type of (categorical) algebra for the $E_2$ operad.

Each MTC $\cC$ determines a projective representation of $\operatorname{SL}(2,\bZ)$ called the \define{modular data} of $\cC$. We will refer to this representation as $\int_{T^2}\cC$; the notation follows the ``factorization algebra'' notation used for example in \cite{WalkerTQFTs,ScheimbauerThesis}. The representation $\int_{T^2}\cC$ has many nice properties. Among them are (i) 
$\int_{T^2}\cC$ has a basis indexed by the simple objects in $\cC$ and (ii) the representation $\int_{T^2}\cC$ is defined over the cyclotomic field $\bQ^{\mathrm{cyc}}$.
Our Main Theorem strongly suggests that, in the case when $\cC = \Rep(W)$ for a regular VOA $W$,

\mytheorem[Conjecture]{At the level of modular data, we have $\int_{T^2} \Rep(^\gamma W) \cong \gamma^2\bigl( \int_{T^2} \Rep(W)\bigr)$.}

This conjecture is enticingly close to the following result of \cite{MR1266785,MR3435813}. Let $\cC$ be an arbitrary MTC, and choose $\gamma \in \Gal(\bQ^{\mathrm{cyc}}/\bQ) \cong \widehat{\bZ}^\times$. Then there is a permutation $\sigma$ of the simple objects of $\cC$, and hence of the basis of $\int_{T^2}\cC$, which intertwines $\int_{T^2}\cC$ with $\gamma^2(\int_{T^2}\cC)$.
Thus the conjecture predicts that $\int_{T^2}\Rep(W)$ and $\int_{T^2}\Rep(^\gamma W)$ are isomorphic as projective $\mathrm{SL}(2,\bZ)$-modules.

\subsection{Cyclotomicity of MTCs}

Suppose one has some arbitrary module $X$ of $\Gal(\bar{\bQ}/\bQ)$ such that the formula $\gamma \triangleright x := \gamma^2(x)$ also defines an action of $\Gal(\bar{\bQ}/\bQ)$ on $X$. Then the action on $X$ factors through $\Gal(\bQ^{\mathrm{cyc}}/\bQ)$. Indeed, if $(\gamma_1\gamma_2)^2 = \gamma_1^2\gamma_2^2$ then $\gamma_1\gamma_2=\gamma_2\gamma_1$, and the abelianization of $\Gal(\bar{\bQ}/\bQ)$ is precisely $\Gal(\bQ^{\mathrm{cyc}}/\bQ)$.

The repeated occurrences of a squared Galois action suggest that perhaps the entire theory of MTCs is cyclotomic. For comparison, finite group theory is cyclotomic --- all characters are cyclotomic, and all representations split over $\bQ^{\mathrm{cyc}}$ --- and MTCs are a version of ``finite quantum groups.'' The following conjectures have been contemplated by various authors, e.g.\ \cite{MR2922607,EvansGannon}.

\mytheorem[Conjectures]{\begin{enumerate}
  \item If $\cC$ is an MTC, then $\cC$ has a (canonical?) form over $\bQ^{\mathrm{cyc}}$.
  \item If $\cC$ is an MTC, then $\cC$ arises as $\Rep(W)$ for some (canonical?) regular VOA $W$.
  \item If $W$ is a regular VOA, then $\Rep(W)$ has a (canonical?) form over $\bQ^{\mathrm{cyc}}$.
\end{enumerate} \vspace*{-9pt}}

\end{document}